\documentclass[a4paper]{article}
\usepackage{hyperref,color,ifthen}
\usepackage{amsmath,amsfonts,amssymb,epsfig}
\usepackage{graphicx}

\newtheorem{prop}{Proposition}[section]
\newtheorem{example}{Example}[section]

\newenvironment{fig}[4]{\begin{center}\includegraphics[width=#3pt,height=#4pt]{#2.#1}}{\end{center}}

\def\pFIG#1,#2,#3,#4,#5{\begin{figure}\begin{fig}{#1}{#2}{#3}{#4}\end{fig}\caption{}\label{#2}\end{figure}}

\def\eFIG#1#2{\begin{figure}\label{#1}\begin{center}\epsfig{file=#1.eps, height=#2cm}\caption{}\end{center}\end{figure}}

\def\set#1{\{#1\}}
\newcommand{\pd}[2][ ]{\ifthenelse{\equal{#1}{ }}{\frac{\partial}{\partial#2}}{\frac{\partial#1}{\partial#2}}}
\def\tseq[#1][#2]#3,#4{#3,#4}
\def\seqp#1,#2,#3{#1_{#2}\dots#1_{#3}}
\def\seqs#1,#2,#3,#4{#1_{#4_{#2}}\dots#1_{#4_{#3}}}

\newcommand{\seq}[2][ ]{\ifthenelse{\equal{#1}{ }}{\seqp#2}{\seqs#2,#1}}
\def\seqmpp#1,#2,#3{#1_{#2}\cdot\dots\cdot#1_{#3}}
\def\seqmps#1,#2,#3,#4{#1_{#4_{#2}}\cdot\dots\cdot#1_{#4_{#3}}}
\def\seqmsp#1,#2,#3,#4{{#1}_{#2}#4\dots#4{#1}_{#3}}
\def\seqmss#1,#2,#3,#4,#5{#1_{#4_{#2}}#5\dots#5#1_{#4_{#3}}}
\newcommand{\seqmp}[2][ ]{\ifthenelse{\equal{#1}{ }}{\seqmpp#2}{\seqmps#2,#1}}
\newcommand{\seqms}[2][ ]{\ifthenelse{\equal{#1}{ }}{\seqmsp#2}{\seqmss#2,#1}}
\newcommand{\seqm}[2][\cdot]{\ifthenelse{\equal{#1}{\cdot}}{\seqmp#2}{\seqms#2,#1}}

\title{Algebraic structure of the space of homotopy classes of cycles and singular homology
\footnote{This work was partially supported by RFBR 01-01-00546 grant.}}
\author{Valery DOLOTIN\footnote{e-mail:\ {\it vd@gate.itep.ru}}\\
Institute of Theoretical \& Experimental Physics, Moscow}

\begin{document}
\maketitle
\vspace{-9cm}
\begin{flushright}
Preprint ITEP-TH-78/02
\end{flushright}
\vspace{7.5cm}
\begin{abstract}
Described the algebraic structure on the space of homotopy classes of cycles with marked topological flags of disks. This space is a non-commutative monoid, with an Abelian quotient corresponding to the group of singular homologies $H_k(M)$. For the marked flag contracted to a point the multiplication becomes commutative and the subgroup of spherical cycles corresponds to the usual homotopy group $\pi_k(M)$.
\end{abstract}
\section{Basic example}
Take $M=\mathbb{R}^2-\{0\}$. Then a non-contractible cycle $c$ in $M$ corresponding to a map of $S^1$ is a submanifold, homeomorphic to $S^1$. For a flag $D$ of two disks ${\rm pt}=D^0\subset D^1$ with $D^0\subset \partial D^1$ we have a subset $[c]_D$ of cycles with $D\subset c$. This subset is a disjoint union of two components $[c]_D=[c]^+_D\cup [c]^-_D$ which differ by the isotopy class of $D\hookrightarrow [c]$ (see Fig.\ref{circ}).

Choosing for each homotopy class $[c]_D$ of flagged cycles $D\hookrightarrow S^1\subset\mathbb{R}^2-\{0\}$ a positive (and complementary negative) component we get a composition law (see Fig.\ref{cmult})
$$[c_1]\circ[c_2]:=[c_1]^+\#_D[c_2]^-$$
\begin{prop}
This composition law defines on the set of homotopy classes of $S^1$-cycles the structure of the group, isomorphic to $\pi_1(\mathbb{R}^2-\{0\})$.
\end{prop}

\eFIG{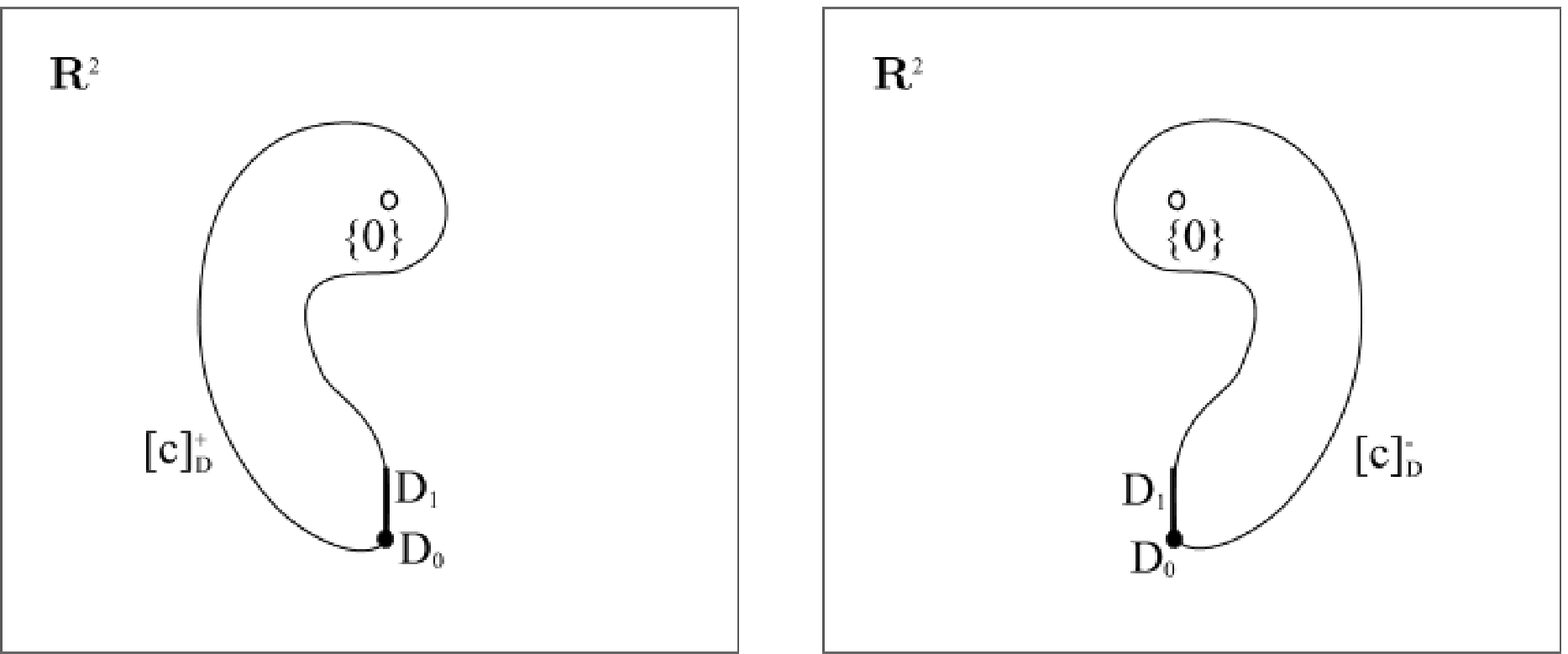}{5}

\eFIG{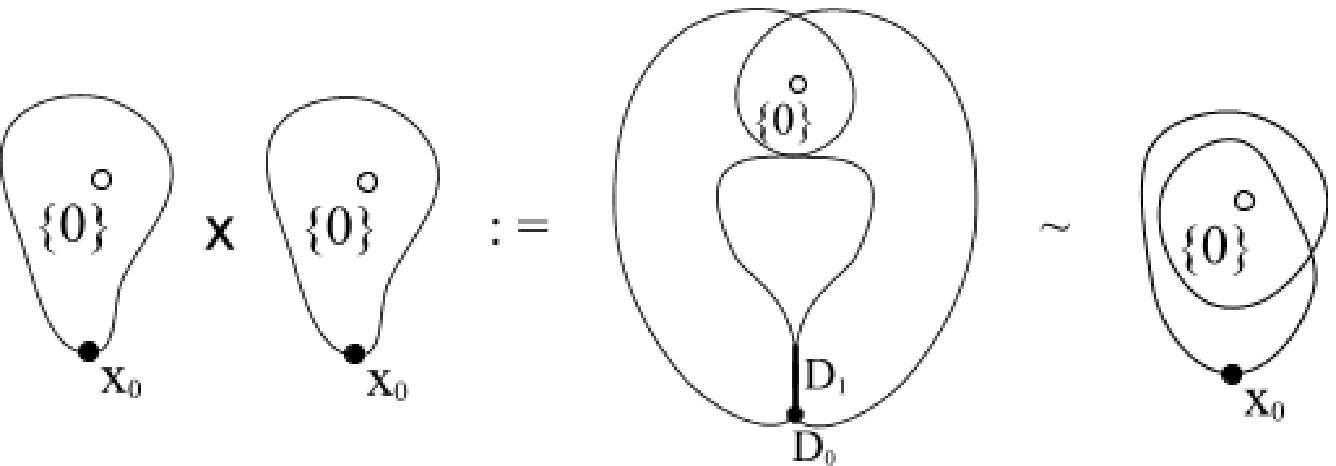}{4}

\section{General construction}

\subsection{Composition law on homotopy classes of cycles}
Fix a flag $D$ of disks $D^0\subset D^1\subset\dots\subset D^k\subset M^n$, where $D^i\subset\partial D^{i+1}$. Denote by $D'$ the subflag $D^0\subset D^1\subset\dots\subset D^{k-1}$ of codimension $1$. By a cycle in $M$ will be called a map $K\to M$ of a closed complex $K$ modulo automorphisms of the domain of the map. Then the homotopy of two cycles is weaker then the homotopy of maps (see Section 1).
Each homotopy class $[c]$ of $k$-cycles on $M$ has a subset of representatives $[c]_{D^i}=\set{c\in [c]\ |\ D^i\subset c}$, i.e. cycles containing $D^i$. In particular $[c]_{D}\subset [c]_{D'}\subset\dots\subset[c]_{D^{(k)}}\subset[c]$. 

Denote by $\Pi_k(M;X)$ the space of homotopy classes of $k$-cycles with a fixed subset $X$.

Fix an orientation $o$ on $D\subset M$. Any two representatives of $[c]_{D'}$ are connected by a $D'$-fixed homotopy, so $[c]_{D'}$ is an element (a class) in $\Pi_k(M;D')$. Take a representative $c\in[c]_{D'}$. Then for a pair of representatives $c_1,c_2\in[c]_{D}\subset[c]_{D'}$ we have a $D'$-fixed homotopy $c_1\stackrel{D'}{\sim} c_2=c_1\stackrel{D'}{\sim} c\stackrel{D'}{\sim} c_2$. This induces an orientation on $D_c\subset c$ (where $D_c$ is $D$ as a submanifold in $c$) for all $c\in[c]_D$. But if for some ("positive") cycle $c^+$ we choose an orientation on $D_{c^+}$ to be equal to $o$ then there is a free (with no fixed set) homotopy of $c^+$ to a cycle $c^-\in[c]_D$ with the homotopy-induced orientation of $D_{c^-}$ opposite to $o$. So the space $[c]_D$  has two components $[c]_D=[c]^+_D\cup[c]^-_D$, where $[c]^+_D$ is the set of cycles with the above homotopy-induced orienatation on $D_c$ equal to $o$, and $[c]^-_D$ - those with the opposite to $o$. Elements of each component $[c]^\sigma_D$ are connected by homotopies with fixed $D$ (i.e. belong to the same class of $\Pi_k(M;D)$), while homotopies between elements from opposite components belong to a wider class of $D'$-fixed homotopy. Then
\begin{prop}
There is a well defined double-covering $s:\Pi_k(M;D)\to\Pi_k(M;D')$, with $s^{-1}([c])=\set{[c]^+_D,[c]^-_D}\subset\Pi_k(M;D)$.
\end{prop}

Choosing a "positive" representative as above fix a splitting $[c]_D=[c]^+_D\cup[c]^-_D$ for each class $[c]\supset[c]_D$, which may be formulated as choosing a "positive section" of the above covering $\Pi_k(M;D)\to\Pi_k(M;D')$.

For each choice of a pair of representatives $a^{\sigma}\in[a]^{\sigma}_D,b^{\tau}\in[b]^{\tau}_D$ we may take their connected sum along the common $D$, denoted by $a^{\sigma}\#^{\tau}b$. After making the pairing of representatives $a^{\sigma}\# b^{\tau}$, we take a class of the result with respect to homotopies with $D'$ fixed. The $D'$-fixed homotopy class of the result does not depend on the choice of representatives incide $[a]^\sigma_D$ and $[b]^\tau_D$, while $[a]^\sigma_D,[b]^\tau_D$ are homotopy invariants of classes $[a]$ and $[b]$.
\begin{prop}
$^{\sigma}\#^{\tau}$, for $\sigma,\tau=+,-$, is a well-defined pairing on the space $\Pi_k(M;D')$ of classes of homotopy of $k$-cycles with fixed $D'$.
\end{prop}
Then we have $4$ different well-defined ways of pairing for homotopy classes of cycles $[a]$ and $[b]$, corresponding to different choices of $\sigma,\tau$  
$$[a]^+\#^+[b],\ [a]^+\#^-[b],\ [a]^-\#^+[b],\ [a]^-\#^-[b]$$

The results of these pairings are in general 4 different elements of $\Pi_k(M;D')$.

\eFIG{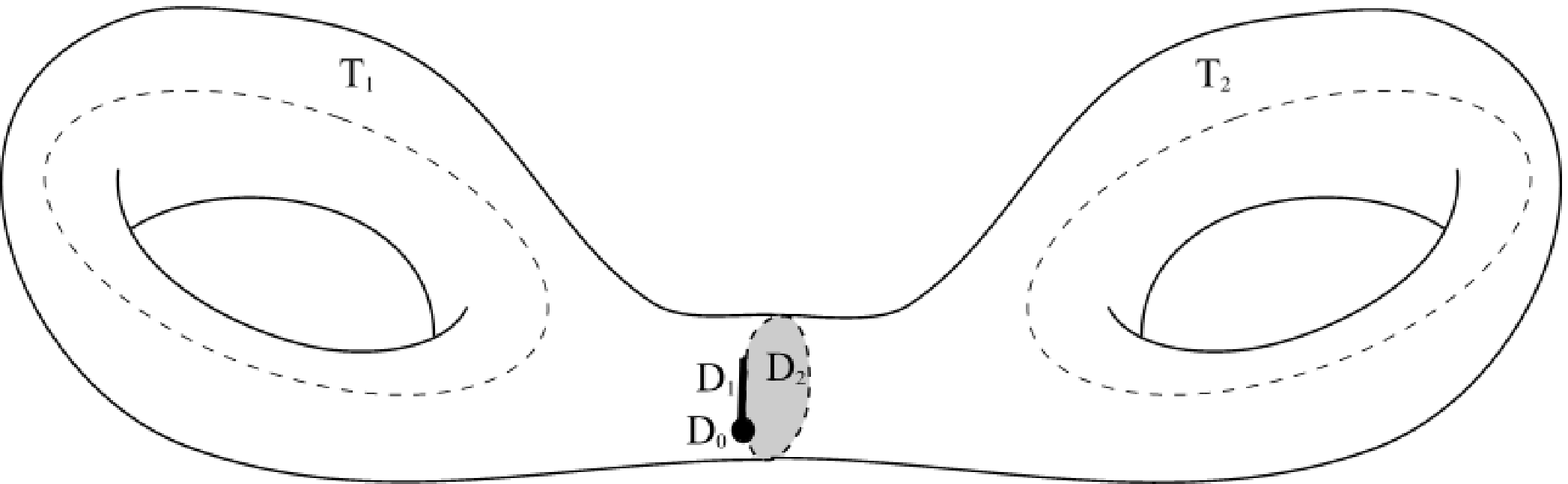}{4}
\eFIG{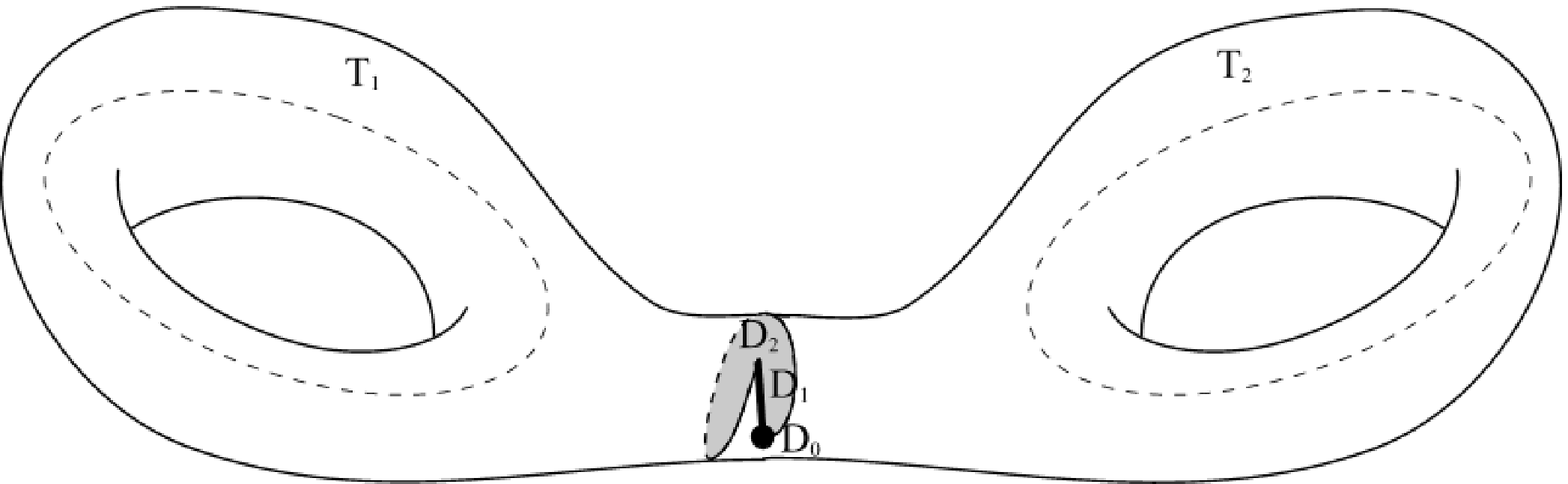}{4}

\begin{example}
Take $M$ to be $\mathbb{R}^3$ with two deleted unlinked circles, $D:D^0\subset D^1\subset D^2,\ D':D^0\subset D^1$, $a:=T_1,\ b:=T_2$ being two tori. Choose decompositions $[T_1]_D=[T_1]^+_D\cup[T_1]^-_D$, $[T_2]_D=[T_2]^+_D\cup[T_2]^-_D$ in a way that the pairing $T_1^+\#^-T_2$ is $D'$-homotopic to that shown on Fig.\ref{toriPM}.

Then the pairing $T_1^-\#^+T_2$ is $D'$-homotopic to that shown on Fig.\ref{toriMP} (here  some $D'$-homotopy has sort of "180-degrees twisted" $D^2$ arownd horizontal axis while $D^1$ being fixed), which has no $D'$-fixed homotopy to $T_1^+\#^-T_2$.

So $[T_1]^+\#^-[T_2]$ and $[T_1]^-\#^+[T_2]$ are defferent elements of $\Pi_2(M;D')$.
\end{example}

Having $[a]^{\sigma}\#^{\tau}[b]\in \Pi_k(M;D')$ we can pair it in 4 different ways with any other class $[c]\in\Pi_k(M;D')$. These combinations of parings may be represented as a composition of "signed" binary trees Fig.\ref{tree}.

\eFIG{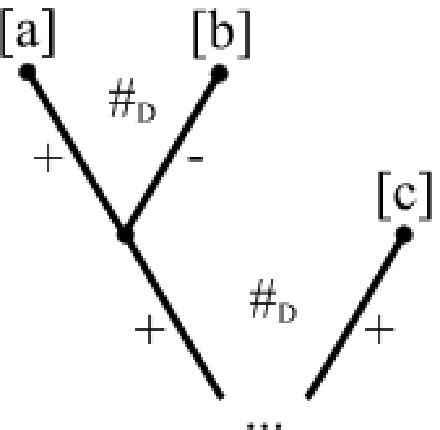}{3}

Note, that $a^{\sigma}\#^{\tau}b\equiv b^{\tau}\#^{\sigma}a$. Then each of the four parings has two equivalent presentations in terms of words of letters $a^+, a^-, b^+, b^-$ as follows:
$$\set{a^+\#^-b\equiv b^-\#^+a}=:\set{ab\sim b^{-}a^{-}}$$
$$\set{a^-\#^+b\equiv b^+\#^-a}=:\set{a^-b^-\sim ba}$$
$$\set{a^+\#^+b\equiv b^+\#^+a}=:\set{ab^-\sim ba^-}$$
$$\set{a^-\#^-b\equiv b^-\#^-a}=:\set{a^-b\sim b^-a}$$
Take an ordered set $\set{[c_1],\dots,[c_N]}$ of generators of $\Pi_k(M;D')$. Let us choose for each pairing of $[c_i],[c_j]$ one presentation (of the two above) in terms of words of letters $c_i^{\sigma_i}$ which will be called "canonical" presentation. For instance, for $i<j$ let $c_ic_j,c_jc_i,c_ic_j^-,c_j^-c_i$ be canonical presentations of the corresponding product classes (connected sums) of $[c_i]$ and $[c_j]$. For a canonical presentation $c=c_{i_1}^{\sigma_1}\dots c_{i_n}^{\sigma_n}$ the complement ("anti-canonical") presentation $c_{i_n}^{-\sigma_n}\dots c_{i_1}^{-\sigma_1}$ will be denoted by $c^-$. If we have the cononical presentations $a$ and $b$ for classes $[a],[b]$, then one (of two) presentatation for $[c]=[a]^{\sigma}\#^{\tau}[b]$ is computed as $a^\sigma b^{-\tau}$.
\begin{example}
\end{example}
Having computed a presentation for the product class we set either it or its complement to be canonical. Then we may use these "signed" presentations of $[c]$ to compute presentations of pairing of $[c]$ with other classes according to the rule above and so on.

{\it Note}, that the ambiguity of choosing one of the two presentations of $[c]$ as the canonical corresponds exactly to the ambiguity while fixing signes of components in the decomposition $[c]_D=[c]_D^+\cup[c]_D^-$.

\begin{prop}\label{TrIso}
The above procedure gives an isomorphism between the space of signed binary trees of pairings of $\set{\seqm[,]{\*[c\*],1,N}}$ and the free monoid $<c_1^\pm,\dots,c_N^\pm>$ generated by letters $c^{\sigma}_i,\ \sigma=+,-$.
\end{prop}

Since the space of signed trees covers $\Pi_k(M;D)$ and the composition there corresponds to connected sums of cycles in $\Pi_k(M;D)$ then we have an epimorphism $p:<c_1^\pm,\dots,c_N^\pm>\to\Pi_k(M;D)$, so
\begin{prop}
$\Pi_k(M;D)$ is a quotient of the free monoid $<c_1^\pm,\dots,c_N^\pm>$.
\end{prop}

To proove the associativity holds we note that the equivalence of presentations $(ab)c=a(bc)$ corresponds to the equality $([a]^+\#^-[b])^+\#^-[c]=[a]^+\#^-([b]^+\#^-[c])$ which may be checked directly to be true for connected sums.

{\it Note}, that $cc^-=[c]^+\#^+[c]$ except for the case of fundamental group $\pi_1$ in general does not belong to the kernel of the above epimorphism, so $\Pi_k(M;D)$ is really a monoid (see Example \ref{NoCont}).

\subsection{Morphism to homology group $\Pi_k(M;D)\to H_k(M)$}
Homotopic cycles are homologous, so we have a well-defined map of sets $h:\Pi_k(M;D)\to H_k(M,\mathbb{Z})$. 
\begin{prop}
For any class $[c]\in\Pi_k(M;D)$ the product $[c]^+\#^+[c]$ is a class of cycles with $h$-image homologous to $0\in H_k(M)$.
\end{prop}
\noindent{\it Proof:\ }
Since $c$ is connected, then $\partial D$ is contractible in $c-D$ to a subcycle $c'\subset c$ of codimension $1$. Then $c^+\#^+c$ is contractible in $M$ to the boundary of a tubular neighbourhood of $c'$.  
$\Box$
\begin{example}\label{NoCont}
Take $M=\mathbb{R}^3-S^1$ and $c=T^2$ to be a torus, having this deleted $S^1$ as its axial circle.


Then the product $T^+\#^+T$ is contractible onto a 1-dimensional cycle and thus is homologous to $0\in H_2(M)$.

\end{example}

According to Proposition \ref{TrIso} we may use elements of the monoid $<c_1^\pm,\dots,c_N^\pm>$ generated by letters $c_i^\sigma,\ \sigma=+,-$ to denote elements of $\Pi_k(M;D)$. We have a map of $<c_i^{\pm}>$ to the abelian group $(c_1,\dots,c_N)$, $c_{i_1}^{\sigma_i}\dots c_{i_n}^{\sigma_n}\mapsto n_1\cdot c_1+\dots+ n_N\cdot c_N$ where $n_i$ is the difference of the numbers of $c_i^+$ and $c_i^-$ entering $c_{i_1}^{\sigma_i}\dots c_{i_n}^{\sigma_n}$.

\begin{prop}
The diagram 
$$\begin{array}{ccc}
<c_1^\pm,\dots,c_N^\pm>&\longrightarrow&(c_1,\dots,c_N)\\
p\downarrow&&\downarrow\\
\Pi_k(M;D)&\stackrel{h}{\longrightarrow}&H_k(M)
\end{array}$$
is commutative. So $h:\Pi_k(M;D){\to} H_k(M)$ is an epimorphism of monoids.
\end{prop}

\section{Discussion}
\begin{enumerate}
\item
Since $[c]_{D^{(i)}}\subset [c]_{D^{(i+1)}}$ then the structure of the monoid $\Pi_k(M;D'')$ of homotopies with the fixed subflag of codimension $2$ will be some quotient of the structure of $\Pi_k(M;D')$ and we have a sequence 
$$\Pi_k(M;D')\to\Pi_k(M;D'')\to\dots\to\Pi_k(M;D^0)\hookleftarrow\pi_k(M)$$
\item
In paper \cite{D1} it was developped a calculus of differential forms with coefficients in (matrix) algebra, with the integration defined over the elements of the group of "flagged spheres". In the current framework the integration of a closed form with coefficients in algebra $\frak{g}$ (higher analogue of flat connection froms) over elements of $\Pi_k(M)$ must give a representation of the homotopy monoid $\Pi_k(M)\to\exp\frak{g}$ in the corresponding Lie group.
\end{enumerate}

\end{document}